\documentclass[12pt]{article}

\usepackage{amsmath}
\usepackage{graphicx}
\usepackage{amssymb}
\usepackage[usenames,dvipsnames]{xcolor}

\usepackage[mathscr]{euscript}

\usepackage{amsthm}

\usepackage{tikz-cd}

\usepackage{xcolor}

\usepackage{enumitem}

\definecolor{ggreen}{rgb}{0.0,0.5,0.0}

\makeatletter
\newcommand{\skipitems}[1]{%
  \addtocounter{\@enumctr}{#1}%
}
\makeatother

\newcommand{\Z}{\mathbb{Z}}

\newcommand{\ra}{\rightarrow}

\newcommand{\tsum}{{\textstyle\sum}}

\newcommand{\Pic}{\text{Pic}}

\newcommand{\mE}{\mathcal{E}}
\newcommand{\sE}{\mathscr{E}}
\newcommand{\mF}{\mathcal{F}}

\newcommand{\mL}{\mathcal{L}}
\newcommand{\mM}{\mathcal{M}}
\newcommand{\mO}{\mathcal{O}}
\newcommand{\PP}{\mathbb{P}}

\newcommand{\Arbarello}{1}
\newcommand{\Atiyah}{2}
\newcommand{\Hartshorne}{3}
\newcommand{\Hazewinkel}{4}
\newcommand{\Oda}{5}
\newcommand{\Okonek}{6}

\usepackage{courier}

\begin{document}

\par\noindent {\centering \textbf{STABILIZATION OF DIRECT IMAGES FOR CURVES}\\\text{ }\par}

\par\noindent {\centering \textsc{fedor bogomolov and spencer f. schrandt}\\\text{ }\par}

\par\noindent \textsc{Abstract.}  We discuss phenomena of stabilization for direct images of line bundles over projective curves mapping onto the projective line, for maps of sufficiently big degree.\\%%

\par\noindent {\centering \textsc{1. Introduction}\\\text{ }\par}

\par\noindent We discuss an effect of ``stabilization'' of direct images of line and vector bundles for maps of curves to $\PP^1$ with increasing degree.  For a simple example, the direct image of the trivial line bundle $\mO$ over $\PP^1$ is a direct sum of $\mO$ and copies of $\mO(-1)$ only.  Similar results can be established for any bundle on $\PP^1$, where the direct image depends only on the degree of the map $\PP^1\ra\PP^1$, and not on the map itself.  We use the fact that the cohomology of the direct image is the same as that of the original bundle; however, different bundles on the base may have the same cohomology, and hence our formulas capture more subtle information than the cohomology alone.  A somewhat similar but more complicated picture holds for curves of arbitrary genus.  We provide corresponding formulas and analyze the asymptotics for arbitrary genus and degrees of maps to $\PP^1$.  We also discuss duality for direct images and related questions in this context.\\

\par\noindent {\centering \textsc{2. Preliminaries}\\\text{ }\par}

\par\noindent Throughout, we consider surjective morphisms $f:X\ra\PP^1$ of varieties with $X$ a smooth irreducible projective curve of genus $g$ over an algebraically closed field $k$.  We use $n=\deg f$ to denote the degree.\\

\par\noindent Before proceeding, recall the basic theory of vector bundles over $\PP^1$.  We have the trivial bundle $\mO$, and let $\mO(-1)$ denote the tautological line bundle over $\PP^1$, which has dual $\mO(1)$, the Serre twisting sheaf.  Put
$$\mO(j)=\begin{cases}\mO(1)^{\otimes j}\text{ if }j>0\text{;}\\\mO\text{ if }j=0\text{;}\\\mO(-1)^{\otimes|j|}\text{ if }j<0\text{.}\end{cases}$$
Then any line bundle over $\PP^1$ is isomorphic to $\mO(j)$ for some $j\in\Z$.  For bundles of arbitrary rank, we have the following classification:\\

\par\noindent\textbf{Theorem.}  (splitting theorem of Grothendieck) Any vector bundle $\mE$ over $\PP^1$ splits as a direct sum of line bundles $\mO(n_j)$, with the summands determined uniquely up to permutation.\\

\par\noindent Together with the above, the key tool for our computations is the cohomology of the $\mO(j)$.  We have:
$$h^0(\PP^1,\mO(j))=\begin{cases}j+1\text{ if }j\geq0\text{;}\\0\text{ if }j<0\text{;}\end{cases}$$
$$h^1(\PP^1,\mO(j))=\begin{cases}0\text{ if }j\geq0\text{;}\\-j-1\text{ if }j<0\text{.}\end{cases}$$
In the simplest cases this alone is enough to determine the direct image, as we shall see below.\\

\par\noindent We also make frequent use of the projection formula for direct images, which in our case says that for any bundles $\mE$ over $X$ and $\mF$ over $\PP^1$, the following formula holds:
$$f_*(f^*\mF\otimes\mE)=f_*\mE\otimes\mF\text{.}$$

\par\noindent Standard references for these or any other preliminary materials are [\Hartshorne] and [\Arbarello].  For more on the splitting theorem of Grothendieck, see for example [\Okonek, 1.2.4] or, over an arbitrary field, [\Hazewinkel].\\

\par\noindent {\centering \textsc{3. Genus zero}\\\text{ }\par}

\par\noindent In the simplest case $X=\PP^1$, we have a complete picture of the direct images.  In particular, we shall find that the direct image depends only on the degree of the map.\\

\par\noindent As stated, sometimes the cohomology only permits one choice for the direct image.  For example:\\

\par\noindent\textbf{Lemma.}  The direct image of the trivial bundle is given by $f_*\mO=\mO\oplus(n-1)\mO(-1)$.\\
\\
\textit{Proof.}  We know $f_*\mO$ is a direct sum of bundles $\mO(n_j)$ for some integers $n_j$, and there are $n$ summands since the degree of $f$ is $n$.  Since $h^0(f_*\mO)=h^0(\mO)=1$ we must have one of the $n_j$ equal to zero and the rest negative.  And due to $h^1(f_*\mO)=h^1(\mO)=0$, all the negative $n_j$ must equal $-1$.  $\square$.\\

\par\noindent More generally, we require the additional data from the previously mentioned projection formula.\\

\par\noindent\textbf{Lemma.}  For $1\leq i<n$, $f_*\mO(i)=(i+1)\mO\oplus(n-i-1)\mO(-1)$.\\
\\
\textit{Proof.}  Again we can write a direct sum $f_*\mO(i)=\oplus_{j=1}^n\mO(n_j)$.  Looking at the first cohomology we must have $n_j\geq-1$ for each $j$, as otherwise $h^1(f_*\mO(i))>0$ whereas $h^1(\mO(i))=0$.\\
By the projection formula, $f_*\mO(i)\otimes\mO(\ell)=f_*(\mO(i)\otimes f^*\mO(\ell))$ for all $\ell\in\Z$, in particular for $\ell=-\max\{n_j\}$.  We know that $f^*\mO(\ell)=\mO(\ell n)$ since the map is degree $n$, and then we write $\tsum_jh^0(\mO(n_j+\ell))=h^0(\mO(i+\ell n))$.\\
If $\ell<0$, then $i+\ell n<0$ and hence the right hand side of the last equation is zero; however, the left hand side should be positive since $n_j+\ell=0$ when $j$ is such that $n_j$ is at maximum.  Therefore $-1\leq n_j\leq0$ for each $j$.  We conclude the required formula by looking at zeroth cohomology: the bundle $\mO$ must appear as a summand exactly $i+1$ times in order to have $h^0(\oplus_j\mO(n_j))=h^0(\mO(i))=i+1$.  $\square$.\\

\par\noindent In fact this determines the direct images of all bundles over $\PP^1$.\\

\par\noindent\textbf{Proposition.}  We have the formula
$$f_*\mO(qn+i)=(i+1)\mO(q)\oplus(n-i-1)\mO(q-1)$$
for any integer $qn+i$ with $q\in\Z$ and $0\leq i<n$.\\
\\
\textit{Proof.}  Write
$$f_*\mO(qn+i)=f_*(f^*\mO(q)\otimes\mO(i))=f_*\mO(i)\otimes\mO(q)$$
using the projection formula, and then apply the preceding lemmas.  $\square$.\\

\par\noindent Remembering Grothendieck's theorem, we can now immediately compute the direct image of any vector bundle $\mE$ over $\PP^1$ once we know how $\mE$ splits into a sum of line bundles; indeed, if $\mE=\oplus_{j=1}^r\mO(n_j)$, then $f_*\mE=\oplus_{j=1}^rf_*\mO(n_j)$.\\
Therefore we obtain that the direct image with respect to a map $\PP^1\ra\PP^1$ depends only on the degree of said map.\\

\par\noindent {\centering \textsc{4. Genus one}\\\text{ }\par}

\par\noindent All direct images can also be computed in the case that $X$ is an elliptic curve.  This is thanks to Atiyah's classification of bundles over $X$ [\Atiyah].  Write $\sE(r,d)$ for the set of isomorphism classes of indecomposable vector bundles over $X$ with rank $r$ and degree $d$.  Then the elements of $\sE(r,d)$ are parametrized by the points of $X$, which are again identified with the elements of $\sE(1,0)=\Pic^0(X)$.  The notation will be conflated somewhat since for a particular bundle $\mE$ we shall also write $\mE\in\sE(r,d)$ to mean that $\mE$ is a representative of a class in $\sE(r,d)$, i.e. that $\mE$ is indecomposable of rank $r$ and degree $d$ over $X$.\\
\\
The cohomology values $h^j(X,\mE)$ for $j=0,1$, $\mE\in\sE(r,d)$ are listed in [\Oda, Section 2].  We put them below for convenience (note that $\mE_r$ denotes an element of the unique class in $\sE(r,0)$ with $h^0(X,\mE_r)\neq0$, and that $\mE_1$ must be isomorphic to the trivial bundle $K_X=\mO$):
$$h^0(X,\mE)=d\text{, }h^1(X,\mE)=0\text{ if }d>0\text{;}$$
$$h^0(X,\mE)=0\text{, }h^1(X,\mE)=|d|\text{ if }d<0\text{;}$$
$$h^0(X,\mE)=h^1(X,\mE)=0\text{ if }d=0\text{ and }\mE\neq\mE_r\text{;}$$
$$h^0(X,\mE_r)=h^1(X,\mE_r)=1\text{.}$$
These cohomology values alone are enough to address some simple cases.\\

\par\noindent\textbf{Lemma.}  We have the direct image $f_*\mE_r=\mO\oplus(rn-2)\mO(-1)\oplus\mO(-2)$.  And for $\mE_r\neq\mE\in\sE(r,0)$, we have the direct image $f_*\mE=rn\mO(-1)$.\\
\\
\textit{Proof.}  In each case, there is only one choice for the direct image with the correct zeroth and first cohomology.  $\square$.\\

\par\noindent We could similarly derive formulas for the direct images of the bundles in $\sE(r,\pm1)$ using only the cohomology, though these results will also follow from the below.  As before, the projection formula allows us to compute the general case.\\

\par\noindent\textbf{Lemma.}  For $1\leq d<rn$ and a bundle $\mE\in\sE(r,d)$, the direct image is given by $f_*\mE=d\mO\oplus(rn-d)\mO(-1)$.\\
\\
\textit{Proof.}  Write $f_*\mE=\oplus_{j=1}^{rn}\mO(n_j)$, and let $\ell=\max\{n_j\}$.  Then we must have $h^0(f_*\mE\otimes\mO(-\ell))>0$ because $f_*\mE\otimes\mO(-\ell)$ has an $\mO$ direct summand.\\
But $\mE\otimes f^*\mO(-\ell)$ is indecomposable of degree $d-rn\ell$, as $\mE$ is indecomposable and $f^*\mO(-\ell)$ is a line bundle.  So now if $\ell>0$ then $d-rn\ell<0$, which makes $h^0(f_*(\mE\otimes f^*\mO(-\ell)))=h^0(\mE\otimes f^*\mO(-\ell))=0$; this is a contradiction to the projection formula $f_*\mE\otimes\mO(-\ell)=f_*(\mE\otimes f^*\mO(-\ell))$.\\
Therefore each $n_j$ is nonpositive.  Now the cohomology of $\mE$ determines everything.  Since $h^0(X,\mE)=d$, we must have $n_j=0$ for exactly $d$ different values of $j$.  And since $h^1(X,\mE)=0$, all of the negative $n_j$ must be $-1$.  $\square$.\\

\par\noindent The above can be extended to bundles of degree not in $[0,rn)$.\\

\par\noindent\textbf{Proposition.}  Given a bundle $\mE\in\sE(r,d)$, choose $q\in\Z$ such that $0\leq d-qrn<rn$ (i.e. let $q$ be the quotient of the division of $d$ by $rn$).  Then
$$f_*\mE=(d-qrn)\mO(q)\oplus((q+1)rn-d)\mO(q-1)$$
unless $\mE\otimes f^*\mO(-q)=\mE_r$, in which case
$$f_*\mE=\mO(q)\oplus(rn-2)\mO(q-1)\oplus\mO(q-2)$$
instead.\\

\par\noindent\textit{Proof.}  The projection formula and the invertibility of $\mO(-q)$ yield
$$f_*\mE=f_*(\mE\otimes f^*\mO(-q))\otimes\mO(q)\text{.}$$
As long as $d\not\equiv0\pmod{rn}$, the direct image of the bundle $\mE\otimes f^*\mO(-q)$ of degree $d-qrn$ is determined by the last Lemma, and we have
$$f_*\mE=(d-qrn)\mO(q)\oplus((q+1)rn-d)\mO(q-1)\text{.}$$
When $d\equiv0\pmod{rn}$, there are two possibilities for the direct image, which correspond to whether $\mE\otimes f^*\mO(-q)\in\sE(r,0)$ is isomorphic to $\mE_r$.  When $\mE\otimes f^*\mO(-q)=\mE_r$,
$$f_*\mE=\mO(q)\oplus(rn-2)\mO(q-1)\oplus\mO(q-2)\text{;}$$
otherwise we have $f_*\mE=rn\mO(q-1)$, which is consistent with the formula obtained in the $d\not\equiv0\pmod{rn}$ case.  $\square$.\\

\par\noindent So the indecomposable bundles of a given rank $r$ and degree $d$ all have the same direct image, determined only by the degree of the map, with the exception that in the case $d\equiv0\pmod{rn}$, those isomorphic to $\mE_r\otimes f^*\mO(-q)$ have a different direct image.  The direct image of any bundle over $X$, once written as a direct sum of indecomposable bundles, is now also determined.  To compute everything we need only the degree of the map $f$ and the pullback $f^*\mO(1)$ (which determines $f^*\mO(\ell)$ for all $\ell\in\Z$).\\

\par\noindent {\centering \textsc{5. Arbitrary genus}\\\text{ }\par}

\par\noindent In low genus we are fortunate enough to be able to compute everything, thanks to nice classifications of the vector bundles over $X$ and known cohomology.  For curves of higher genus, we cannot expect such a simple and definitive description of the direct images.  However, we still observe an effect of ``stabilization,'' wherein the degrees of the summands $\mO(n_j)$ which appear in the direct image become increasingly concentrated as the degree of the map becomes large.  We have already seen that when $g=0$ or $1$, the direct image of an indecomposable bundle has summands of degree $q$ and $q-1$ only for some $q\in\Z$, with one special situation in the $g=1$ case where a single summand of degree $q-2$ is present (this is even irrespective of the degree of the map).  For $g>1$ the phenomenon is somewhat more complicated, and even more so for higher rank bundles, but focusing on line bundles we can analyze the asymptotics of this effect.  Prior to this discussion, we make a detour to discuss the dualizing behavior of direct images; for this we can consider arbitrary bundles over $X$, and the result shall be useful in that which follows it.\\

\par\noindent ``Duality'' in our case is the phenomenon that the direct images themselves mirror the dualizing feature of the cohomology given by Serre.  That is, the action of taking the dual and then tensoring with the canonical bundle commutes with taking the direct image.  Indeed, Serre duality itself suggests the result since it implies that the two relevant bundles have the same cohomology.  We may also confirm via our previous work that it holds already for $g=0$ or $1$.  Another simple example for arbitrary genus is the case where $n=2$ and $\mL$ is a line bundle with $h^0(\mL),h^1(\mL)\neq0$.  There, $f_*\mL=\mO(h^0(\mL)-1)\oplus\mO(-h^1(\mL)-1)$, and since $\mM=K_X\otimes\mL^*$ has $h^0(\mM)=h^1(\mL)$ and $h^1(\mM)=h^0(\mL)$, we get $f_*\mM=\mO(h^1(\mL)-1)\oplus\mO(-h^0(\mL)-1)$.  So we see that $f_*\mM=K_{\PP^1}\otimes(f_*\mL)^*$.\\

\par\noindent In the course of the proof below, we shall write the coefficients $n_j$ in the direct image $f_*\mE=\oplus_{j\in\Z}n_j\mO(j)$ of a bundle $\mE$ in terms of the cohomology of the bundles $f^*\mO(\ell)\otimes\mE$, whose direct images are the twists $f_*\mE\otimes\mO(\ell)$ of the original direct image.  We motivate our strategy as follows.\\

\par\noindent By twisting $f_*\mE$ to a sufficient degree, we can ensure that the highest-degree summands are copies of $\mO$, so that the zeroth cohomology of the twisted bundle tells us the number of those summands, hence the number of highest-degree summands in $f_*\mE$.  Then, by twisting one step less so that the highest-degree summands are copies of $\mO(1)$, the zeroth cohomology of the resulting bundle will allow us to again count the number of $\mO$ summands, which is the number of summands of second-highest degree.  Continuing in this manner, we compute all the $n_j$ in terms of the zeroth cohomology of the twisted bundles.  Shifting the other way, we could alternatively write the $n_j$ in terms of the first cohomology of the twisted bundles: we would shift so that the lowest-degree summands not yet counted are copies of $\mO(-2)$.  Doing this both for $\mE$ and for $K_X\otimes\mE^*$ allows us to compare the coefficients and obtain the result.\\

\par\noindent\textbf{Theorem.}  Let $f:X\ra\PP^1$ be a degree $n$ surjective map of smooth algebraic curves.  Then
$$f_*(K_X\otimes\mE^*)=K_{\PP^1}\otimes(f_*\mE)^*$$
for any vector bundle $\mE$ over $X$.\\
\\
\textit{Proof.}  Let $\mE$ have rank $r$.  We know there are splittings
$$f_*\mE=\oplus_{j\in\Z}n_j\mO(j)\text{, }f_*(K\otimes\mE^*)=\oplus_{j\in\Z}m_j\mO(j)$$
with $n_j,m_j\geq0$ and $\tsum_jn_j=\tsum_jm_j=rn$.\\
Since the cohomology groups are finite-dimensional, there are only finitely many $j\in\Z$ for which the $n_j,m_j$ can possibly be nonzero.  Set $M=h^0(\mE)-1$ and $M'=-h^1(\mE)-1$; also note that $-M-2=-h^1(K\otimes\mE^*)-1$ and $-M'-2=h^0(K\otimes\mE^*)-1$.  Then we have $n_j=0$ if $j>M$ or $j<M'$; and $m_j=0$ if $j>-M'-2$ or $j<-M-2$.\\
Now define the numbers $a_\ell=h^0(f^*\mO(-\ell)\otimes\mE)$ for $\ell\in\Z$.  Since
$$f_*(f^*\mO(-\ell)\otimes\mE)=f_*\mE\otimes\mO(-\ell)=\oplus_{j=M'}^Mn_j\mO(j-\ell)\text{,}$$
we have $a_\ell=0$ for $\ell>M$ and
$$a_\ell=\tsum_{j=\ell}^M(j-\ell+1)n_j$$
for $M'\leq\ell\leq M$.  From this we see immediately that $n_M=a_M=a_M-2a_{M+1}+a_{M+2}$.  When $M'\leq\ell<M$ we solve for
$$n_\ell=a_\ell-\tsum_{j=\ell+1}^M(j-\ell+1)n_j\text{;}$$
hence we have written $n_\ell$ in terms of the $n_j$ with $j>\ell$.\\
Now again let $M'\leq\ell<M$ and assume for induction that $n_j=a_j-2a_{j+1}+a_{j+2}$ for all $\ell<j\leq M$.  The equation for $n_\ell$ above can then be rewritten as
$$n_\ell=a_\ell-\tsum_{j=\ell+1}^M(j-\ell+1)a_j+2\tsum_{j=\ell+2}^{M+1}(j-\ell)a_j-\tsum_{j=\ell+3}^{M+2}(j-\ell-1)a_j$$
$$=a_\ell-2a_{\ell+1}+a_{\ell+2}\text{.}$$
By induction we conclude that $n_j=a_j-2a_{j+1}+a_{j+2}$ whenever $M'\leq j\leq M$.\\
To get the required result, we shall also write the $m_j$ in terms of the numbers $a_\ell$; to do this, we perform an essentially dual argument relying on the relation $a_\ell=h^1(f^*\mO(\ell)\otimes K\otimes\mE^*)$.\\
First, reuse the projection formula
$$f_*(f^*\mO(-\ell-2)\otimes K\otimes\mE^*)=f_*(K\otimes\mE^*)\otimes\mO(-\ell-2)$$
$$=\oplus_{j=-M-2}^{-M'-2}m_j\mO(j-\ell-2)\text{.}$$
From this we get again $a_{-\ell-2}=0$ for $\ell<-M-2$, but also the new formula
$$a_{-\ell-2}=\tsum_{j=-M-2}^\ell(\ell-j+1)m_j$$
for $-M-2\leq\ell\leq-M'-2$.  We see in particular that $m_{-M-2}=a_M=a_M-2a_{M+1}+a_{M+2}$.\\
When $-M-2<\ell\leq-M'-2$, write
$$m_\ell=a_{-\ell-2}-\tsum_{j=-M-2}^{\ell-1}(\ell-j+1)m_j$$
and assume for induction that $m_j=a_{-j-2}-2a_{-j-1}+a_{-j}$ for $-M-2\leq j<\ell$.  In this case
$$m_\ell=a_{-\ell-2}-\tsum_{j=-M-2}^{\ell-1}(\ell-j+1)a_{-j-2}$$
$$+2\tsum_{j=-M-3}^{\ell-2}(\ell-j)a_{-j-2}-\tsum_{j=-M-4}^{\ell-3}(\ell-j-1)a_{-j-2}\text{;}$$
hence $m_\ell=a_{-\ell-2}-2a_{-\ell-1}+a_{-\ell}$.  We conclude that $m_j=a_{-j-2}-2a_{-j-1}+a_{-j}$ for $-M-2\leq j\leq-M'-2$, and thus we have written the $m_j$ in terms of the $a_\ell$.\\
The result now follows, since
$$\mO(-2)\otimes(f_*\mE)^*=\oplus_{j=M'}^Mn_j\mO(-j-2)=\oplus_{j=-M-2}^{-M'-2}n_{-j-2}\mO(j)$$
and evidently $n_{-j-2}=m_j$ for $-M-2\leq j\leq-M'-2$.  $\square$.\\

\par\noindent Instead of computing the $m_j$ with the same method as that applied for the $n_j$, we could also use Riemann-Roch to complete the latter part of the proof.  Before, we wrote $n_j=a_j-2a_{j+1}+a_{j+2}$ in terms of the numbers $a_\ell=h^0(f^*\mO(-\ell)\otimes\mE)$; we can switch from zeroth to first cohomology, that is, write $n_j$ in terms of the numbers $b_\ell=h^1(f^*\mO(-\ell)\otimes\mE)$, as follows.  Given $\mE$ with rank $r$ and degree $d$, compute $a_\ell-b_\ell=(d-rn\ell)+r(1-g)$ using Riemann-Roch.  Therefore $(a_j-2a_{j+1}+a_{j+2})-(b_j-2b_{j+1}+b_{j+2})=0$, which implies $n_j=b_j-2b_{j+1}+b_{j+2}$.  Apply this to $K\otimes\mE^*$ and we recover our formula for $m_j$ in terms of the $a_\ell$.\\

\par\noindent We now address the previously posed question of stability.  For this, restrict to the case of a line bundle $\mL$ over $X$ with $d=\deg\mL$, and choose $q\in\Z$ such that $0\leq d-qn<n$.  We may as well assume $g>1$ from here on out.  Given the decomposition $f_*\mL=\oplus_{j=1}^n\mO(n_j)$, we shall like to discuss the range of the degrees of the summands, and hence define $s(\mL,f)=\max_{1\leq j,j'\leq n}|n_j-n_{j'}|$.\\

\par\noindent We will say the direct images become ``stable'' in a sense that, when the degree of the map grows big, the degrees of the summands $\mO(n_j)$ which can appear are limited to a small range, i.e. $s(\mL,f)$ has a small bound.  In the most extreme case, taking $n$ excessively large, this effect is the following: for $0\leq d<n$, the summands which account for the cohomology must be as numerous as possible, that is split into $h^0(\mL)$ copies of $\mO$ and $h^1(\mL)$ copies of $\mO(-2)$, while the rest and vast majority of the summands are copies of $\mO(-1)$ which contribute nothing to the cohomology; due to the projection formula, for arbitrary degrees $d$ of the line bundle the same occurs except with the degrees of all summands ``shifted'' by a fixed amount, exactly the consequence of tensoring the direct sum by some $\mO(\ell)$.  This example is expressed by the following result, after which we shall see similar phenomena even as the condition on the size of the degree of the map is eased.\\

\par\noindent\textbf{Proposition.}  If $n>2g-2$, then $s(\mL,f)\leq2$.\\
\\
\textit{Proof.}  First assume $0\leq d<n$.  Then since $f^*\mO(-1)\otimes\mL$ has degree $d-n<0$, $f_*\mL\otimes\mO(-1)=f_*(f^*\mO(-1)\otimes\mL)$ has no nonzero global sections; and since $f^*\mO(1)\otimes\mL$ has degree $d+n>2g-2$, $f_*\mL\otimes\mO(1)=f_*(f^*\mO(1)\otimes\mL)$ is non-special.  We conclude that $f_*\mL$ cannot have any summands of positive degree or of degree less than $-2$.  This means the direct image must be
$$f_*\mL=h^0(\mL)\mO\oplus(n-h^0(\mL)-h^1(\mL))\mO(-1)\oplus h^1(\mL)\mO(-2)$$
as described above.\\
For arbitrary $d$, we apply the result to $f^*\mO(-q)\otimes\mL$ which has degree $d-qn$; then $f_*\mL=f_*(f^*\mO(-q)\otimes\mL)\otimes\mO(q)$.  $\square$.\\

\par\noindent\textbf{Corollary.}  If $n>2$, $g=2$, then $s(\mL,f)\leq2$.\\
\\
\textit{Proof.}  If $g=2$ then $2g-2=2$.  $\square$.\\

\par\noindent We shall address the case of $n=2$ for general $g>1$ below.  But first, we discuss the weaker condition $n>g-1$.  For general bundles, even this is enough to get the best possible.\\

\par\noindent\textbf{Proposition.}  Say $n>g-1$ and that $\mL$ is a \textit{general} line bundle.  Then $s(\mL,f)\leq1$.\\
\\
\textit{Proof.}  For a general line bundle $\mM$ of degree $\deg\mM\geq g-1$, we have $h^1(\mM)=0$.  And, dually, for general $\mM$ of degree $\deg\mM\leq g-1$ we have $h^0(\mM)=0$.\\
So let $\mL$ be general with $g-1\leq d<g-1+n$.  Then $h^1(\mL)=0$ and so all summands of $f_*\mL$ must be of degree at least $-1$.  Furthermore, $f_*\mL$ has no summands of positive degree since $f^*\mO(-1)\otimes\mL$ has no nonzero global sections (meaning the same is true for $f_*\mL\otimes \mO(-1)$).  This leaves only the possibility
$$f_*\mL=h^0(\mL)\mO\oplus(n-h^0(\mL))\mO(-1)\text{.}$$
So $s(\mL,f)\leq1$ holds for general $\mL$ with $g-1\leq d<g-1+n$.  Therefore the bound holds for all $d$ since for general $\mL$ of arbitrary degree we can apply the result for some $f^*\mO(\ell)\otimes\mL$.  $\square$.\\

\par\noindent Removing the requirement that $\mL$ be general, we manage some partial results.  Relying again on the projection formula and conditions under which the cohomology vanishes, we see that the best chance of results comes for bundles of degree close to $g-1$, where said degree is not too far from zero or from the value $2g-2$.\\

\par\noindent\textbf{Proposition.}  Suppose $n>g-1$.  Then if $d=g-1$, then $s(\mL,f)\leq2$.  If $d=g-2$, then $s(\mL,f)\leq3$ with equality only possible when $\mL=f^*\mO(-1)\otimes K$ (in which case $n=g$).  If $d=g$, then $s(\mL,f)\leq3$ with equality only possible when $\mL=f^*\mO(1)$ (in which case $n=g$ once again).\\
\\
\textit{Proof.}  If $d=g-1$, then since $d+n>0$ and $d-n<0$, the direct image cannot have any summands of degree greater than $0$ or less than $-2$.  So again
$$f_*\mL=h^0(\mL)\mO\oplus(n-h^0(\mL)-h^1(\mL))\mO(-1)\oplus h^1(\mL)\mO(-2)$$
and $s(\mL,f)\leq2$.\\
If $d=g-2$, then $d-n<0$ and so the direct image of $\mL$ cannot have any summands of degree greater than $0$.  Also $d+n\geq2g-2$, so $h^1(f_*\mL\otimes\mO(1))=1$ or $0$, depending on if $f^*\mO(1)\otimes\mL=K$ or not.  In the first case we must have
$$f_*\mL=h^0(\mL)\mO\oplus(n-h^0(\mL)-h^1(\mL)+1)\mO(-1)\oplus(h^1(\mL)-2)\mO(-2)\oplus\mO(-3)$$
and in the second case we must have
$$f_*\mL=h^0(\mL)\mO\oplus(n-h^0(\mL)-h^1(\mL))\mO(-1)\oplus h^1(\mL)\mO(-2)\text{.}$$
In any case $s(\mL,f)\leq3$, with equality only possible for $f^*\mO(1)\otimes\mL=K$.\\
The case $d=g$ is exactly dual to the case of $d=g-2$.  Given $\mL$ of degree $g$, form $K\otimes\mL^*$ of degree $g-2$ and use the above to calculate its direct image; then finish with $f_*(K\otimes\mL^*)=\mO(-2)\otimes(f_*\mL)^*$, which allows us to compute $f_*\mL$.  Note also that the condition for equality in the bound on $s(\mL,f)$ in this case is $K\otimes\mL^*=f^*\mO(-1)\otimes K$, i.e. $\mL=f^*\mO(1)$.  $\square$.\\

\par\noindent Lastly, we apply our established methods to arbitrary $n$.  Afterwards, we analyze the asymptotics for large $g,n$.\\

\par\noindent\textbf{Proposition.}  If $n$ is odd then $s(\mL,f)\leq\frac{2g-5/2}{n}+\frac{5}{2}$; or generally, $s(\mL,f)\leq\frac{2}{3}g+\frac{5}{3}$.  If $n$ is even then $s(\mL,f)\leq\frac{2g-2}{n}+\frac{5}{2}$.  For $n=2$ this gives $s(\mL,f)\leq g+1$, and for $n\neq2$ even we get $s(\mL,f)\leq\frac{1}{2}g+2$.  Therefore if $X$ is not hyperelliptic to begin with or if $f$ is not the hyperelliptic map $X\ra\PP^1$ of degree $2$, then the bound $s(\mL,f)\leq\frac{2}{3}g+\frac{5}{3}$ holds for all $n$.\\
\\
\textit{Proof.}  First let $n$ be odd.  Then it suffices to show the result for $g-1-\tfrac{n-1}{2}\leq d\leq g-1$, for then duality shall give us the same bound for $g-1\leq d\leq g-1+\tfrac{n-1}{2}$, and we shall have it for $n$ consecutive values of $d$.\\
We note that if $\ell>\tfrac{g-3/2}{n}+\tfrac{1}{2}$, then $(g-1-\tfrac{n-1}{2})+\ell n>2g-2$.  Therefore no $n_j$ can be strictly less than
$$-\left(\frac{g-3/2}{n}+\frac{1}{2}+1\right)-1=-\left(\frac{g-3/2}{n}+\frac{5}{2}\right)\text{;}$$
for otherwise $f_*\mL\otimes\mO(\ell)$ would have a summand of degree less than $-1$, where $\ell=\left[\frac{g-3/2}{n}+\frac{1}{2}\right]+1$.  On the other side, if $\ell>\tfrac{g-1}{n}$, then $(g-1)-\ell n<0$.  So no $n_j$ can be strictly bigger than
$$\left(\frac{g-1}{n}+1\right)-1=\frac{g-1}{n}\text{;}$$
for otherwise $f_*\mL\otimes\mO(-\ell)$ would have a summand of nonnegative degree, where $\ell=\left[\tfrac{g-1}{n}\right]+1$.  Together, we have
$$s(\mL,f)\leq\frac{g-1}{n}+\left(\frac{g-3/2}{n}+\frac{5}{2}\right)=\frac{2g-5/2}{n}+\frac{5}{2}\text{.}$$
Since $n\geq3$, the right hand side above is less than or equal to $\tfrac{2}{3}g+\tfrac{5}{3}$.\\
So now let $n$ be even, and we perform a similar argument.  It suffices to show the result for $g-1-\tfrac{n}{2}\leq d\leq g-1$, and then duality gives the same bound for $g-1\leq d\leq g-1+\tfrac{n}{2}$.\\
We note that if $\ell>\tfrac{g-1}{n}+\tfrac{1}{2}$, then $(g-1-\tfrac{n}{2})+\ell n>2g-2$.  Therefore no $n_j$ can be strictly less than
$$-\left(\frac{g-1}{n}+\frac{1}{2}+1\right)-1=-\left(\frac{g-1}{n}+\frac{5}{2}\right)\text{.}$$
And just as before, if $\ell>\tfrac{g-1}{n}$, then $(g-1)-\ell n<0$, so no $n_j$ can be strictly greater than
$$\left(\frac{g-1}{n}+1\right)-1=\frac{g-1}{n}\text{.}$$
We conclude that
$$s(\mL,f)\leq\frac{g-1}{n}+\left(\frac{g-1}{n}+\frac{5}{2}\right)=\frac{2g-2}{n}+\frac{5}{2}\text{.}$$
If $n\neq2$, then $n\geq4$ so replacing $n$ with $4$ above we get $s(\mL,f)\leq\tfrac{1}{2}g+2$.  $\square$.\\

\par\noindent Of course, we have already seen that when $n>2g-2$, there is the bound $s(\mL,f)\leq2$; this is optimal since a bundle $\mL$ with both $h^0(\mL),h^1(\mL)$ nonzero will have $s(\mL,f)\geq2$.  Otherwise, assuming $g,n$ are big, the above tells us that the range $s(\mL,f)$ of degrees of the summands which can appear exceeds the ideal bound of $2$ by roughly the quotient of $2g$ by $n$.  \\

\par\noindent Thoughout, we have seen that the degrees of the summands which appear in the direct image are forced to be very ``small,'' that is, relative to the quotient of $d$ by $n$.  In a sense, then, these direct images are approximate to the most ``generic'' bundle on $\PP^1$ with the same cohomology.  In the case of $n>2g-2$ sufficiently large, or any other with $s(\mL,f)\leq2$, this statement is precise: since the direct image is just a shifted version of the case where all summands $\mO(n_j)$ have $-2\leq n_j\leq0$, for some $q\in\Z$ in the general case we'll have $\mO(q)$ and $\mO(q-2)$ summands which account for the cohomology of the shifted bundle, with the rest of the summands simply $\mO(q-1)$.  But even when we do not have $s(\mL,f)\leq2$, we still receive a direct image which is roughly similar.  Recall also that the cohomology of the original bundle does not alone determine the direct image except in some simple cases, and moreover very little beyond the degree of the map comes into play, so our computations capture new, potentially quite subtle information about the bundle.  We think these phenomena should extend to other nontrivial situations.\\

\par\noindent {\centering \textsc{References}\\\text{ }\par}

\begin{enumerate}[label={[\arabic*]},leftmargin=*]

\item E. Arbarello, M. Cornalba, P. A. Griffiths, and J. Harris, Geometry of algebraic curves volume I, Grundlehren der Mathematischen Wissenschaften, 267, Springer, New York, 1985.

\item M. F. Atiyah, Vector bundles over an elliptic curve, Proc. London Math. Soc., \textbf{7} (1957), 414-452.

\item R. Hartshorne, Algebraic geometry, Graduate Texts in Mathematics, 52, Springer, New York, 1977.

\item M. Hazewinkel and C. F. Martin, A short elementary proof of Grothendieck's theorem on algebraic vectorbundles over the projective line, J. Pure and Applied Algebra, \textbf{25} (1982), 207-211.

\item T. Oda, Vector bundles on an elliptic curve, Nagoya Math. J., \textbf{43} (1971), 41-72.

\item C. Okonek, M. Schneider, and H. Spindler, Vector bundles on complex projective spaces, Progress in Mathematics, 3, Birkh$\ddot{\text{a}}$user, Boston, 1980.

\end{enumerate}

\par\noindent \textsc{Fedor A. Bogomolov}\\
\textsc{Department of Mathematics}\\
\textsc{Courant Institute, NYU}\\
\textsc{251 Mercer Street}\\
\textsc{New York, NY 10012, USA,}\\
\texttt{bogomolov@cims.nyu.edu}\textsc{.}\\
\textsc{also:}\\
\textsc{Laboratory of Algebraic Geometry,}\\
\textsc{National Research University Higher School of Economics,}\\
\textsc{Department of Mathematics, 6 Usacheva Street,}\\
\textsc{Moscow, Russia.}\\

\par\noindent \textsc{Spencer F. Schrandt}\\
\textsc{Department of Mathematics}\\
\textsc{Courant Institute, NYU}\\
\textsc{251 Mercer Street}\\
\textsc{New York, NY 10012, USA,}\\
\texttt{sfs9775@nyu.edu}\textsc{.}

\end{document}